\documentclass{rspublic}[]

\usepackage{graphicx}

\newcommand\RR{\mathbb{R}}

\newcommand{\half}{{\mbox{\(\frac{1}{2}\)}}}

\renewcommand{\[}{\begin{equation}}
\renewcommand{\]}{\end{equation}}

\newlength{\figurewidth}
\begin{document}

\title{On Taylor dispersion in oscillatory channel flows}

\author[K. M. Jansons]{Kalvis M. Jansons%
\footnote{Email: \texttt{Dispersion@kalvis.com}}}
\affiliation{Department of Mathematics,
University College London,
Gower Street,
London WC1E 6BT, UK}
\label{firstpage}

\maketitle

\begin{abstract}{Taylor dispersion, oscillatory flows, non-Newtonian
    fluids, power-law fluids, particle sorting, transport in bones and
    connective tissue, zero Reynolds number}{ We revisit Taylor
    dispersion in oscillatory flows at zero Reynolds number, giving an
    alternative method of calculating the Taylor
    dispersivity that is easier to use with computer algebra packages
    to obtain exact expressions.  We consider the effect of
    out-of-phase oscillatory shear and Poiseuille flow, and show that
    the resulting Taylor dispersivity is independent of the phase
    difference.  We also determine exact expressions for several
    examples of oscillatory power-law fluid flows.}
\end{abstract}

\section{Introduction}

The standard Taylor-dispersion problem (Taylor (1953)) is to determine
the approximate long-time distribution of some passive tracer with
density $\phi$ in a channel or pipe $\Omega = \RR\times S$, where $S$
is some connected smooth one- or two-dimensional manifold.  The
governing equation in PDE form is
\[\label{general}
\frac{\partial\phi}{\partial t} + \nabla\cdot(V\phi - D\nabla\phi) = 0,
\]
where $t$ is time, $V$ is the velocity field, and $D$ is the
diffusivity, usually with zero-flux boundary conditions on the walls.
Here we focus on the case where $S=[0,a]$ and the component of the
velocity along the channel is a periodic function of time.

The early work on oscillatory flows focused on flow in pipes (Aris
(1960)), and much later work on oceanographic applications in
unbounded flows (Young et. al. (1982)), but more recently there has
been interest in biomechanics, notably Schmidt et. al. (2005), and the
references therein, who studied consequences of oscillatory flows on
nutrient transport in bone during physical activity and ultrasonic
therapy.  In these biomechanical applications, the flow is probably
better modelled by a combination of shear- and pressure-driven flows,
and natural questions are whether these contributions interact, and
whether the phase difference between the two components is important.

Also in this study we find an expression for the Taylor dispersivity
that is more appropriate for use with computer algebra packages.  This
extends the work of Jansons and Rogers (1995), by providing a more
accessible version of the solution in a form that could be easily
applied by people not familiar with the theory of stochastic
processes, although it uses stochastic calculus in its derivation.

Using computer algebra, we then find the general expression for the
Taylor dispersivity at zero Reynolds number for sinusoidal shear and
Poiseuille flow with arbitrary phase difference, and show in
particular that the Taylor dispersivity is independent of this phase
difference. This extends the work of Claes and Van de Broeck (1991).
We also find exact expressions for the Taylor dispersivity at zero
Reynolds number for several pressure-driven power-law fluid flows.

\section{Dispersion in an oscillatory flow}

We suppose that the \emph{vertical} diffusion $Y$ is confined to the
interval $[0,a]$, which in It\^o form is governed by
\[\label{SDE}
dY_t = \sigma dB_t + v(Y_t)dt + dL^0_t - dL^a_t,
\]
where $B$ is a standard Brownian motion (that is, it has diffusivity
$\half$), and $L^0$ and $L^a$ are local-time terms to keep the
particle in the interval.  (The local times change only when $Y$
hits a boundary.) The above equation describes, in fluid mechanical
terms, a particle with diffusivity $D=\sigma^2/2$ and vertical drift
$v$.  For simplicity, we take $\sigma$ to be constant, as it is in
almost all applications, though this condition could easily be
relaxed, and we also, for simplicity, assume that $v$ is bounded.

As is usual in Taylor-dispersion problems, we neglect the effect of
particle diffusion in the $x$ direction, though it would not be
difficult to include this effect if it were considered useful.  Thus
the \emph{horizontal} displacement is given by
\[
X_t \equiv \int_0^t u(s, Y_s) ds,
\]
where $u(t, y)$ is a periodic function of time, with period
$2\pi/\omega$ say, and ``$\equiv$'' means equal by definition.

From \eqref{SDE}, we see that the stationary vertical density $q$ is
given by
\[
q(y)dy \equiv \lim_{t\to\infty} P[Y_t \in dy] 
= \alpha \exp\left(\int_0^y \frac{2v(\eta)}{\sigma^2} d\eta\right)\,dy,
\]
where $\alpha$ is a normalization constant (so $\int_0^a q(y) dy = 1$).

The part of the horizontal displacement due to the cross-channel mean
of horizontal advection is trivial, and unimportant in the long-time
limit as it is bounded, so we remove this contribution from the
problem by replacing $u(t,y)$ with $u'(t,y) \equiv u(t,y) - \int_0^a
q(y) u(t,y) dy$.

To find the variance of $X_t$, apply It\^o's formula:
\[
df(t,Y_t) = \frac{\partial f}{\partial y} \left(\sigma dB_t + dL^0_t -
dL^a_t\right) + \left(\half\sigma^2 \frac{\partial^2 f}{\partial y^2} +
v \frac{\partial f}{\partial y} + \frac{\partial f}{\partial t}\right) dt,
\]
for some function $f$ to be identified.  Now take $\frac{\partial
  f}{\partial y}(t,0) = \frac{\partial f}{\partial y}(t,a) = 0$ for
all $t$ to remove the local-time terms, and
\[\label{PDE}
\half\sigma^2 \frac{\partial^2 f}{\partial y^2} + v \frac{\partial
  f}{\partial y} + \frac{\partial f}{\partial t} =
u'(t,y),
\]
where $f$ is taken to be an oscillatory function of $t$, with period
$2\pi/\omega$. This gives
\[\label{Xdash}
f(t,Y_t) - f(0,Y_0) = \int_0^t \frac{\partial f}{\partial y}(s,Y_s) dB_s
+ X'_t,
\]
where
\[
X'_t \equiv \int_0^t u'(s, Y_s) ds.
\]
The integral term can be rewritten, namely
\[
\int_0^t \frac{\partial f}{\partial y}(s,Y_s) dB_s = W\left(\int_0^t
\left(\frac{\partial f}{\partial y}(s,Y_s)\right)^2 ds\right),
\]
where $W$ is a standard Brownian motion.  Then the ergodic theorem tell us that 
in the limit $t\to\infty$,
\[
\frac1t \int_0^t
\left(\frac{\partial f}{\partial y}(s,Y_s)\right)^2 ds
\to \frac\omega{2\pi}\int_0^{2\pi/\omega}\int_0^a \left(\frac{\partial
    f}{\partial y}(t,y)\right)^2 q(y) dy \, dt.
\]
First observe that
\[
\lim_{t\to\infty}E[X_t^2]/t = \lim_{t\to\infty}E[X'_t]/t,
\]
and second that $f$ is bounded.  Thus, \eqref{Xdash}, together with
the other observations, implies that the Taylor dispersivity is given
by 
\[
\label{Taylor}
{\cal D} \equiv \lim_{t\to\infty}E[X_t^2]/t =
\frac\omega{2\pi}\int_0^{2\pi/\omega}\int_0^a \left(\frac{\partial
    f}{\partial y}(t,y)\right)^2 q(y) dy \, dt.
\]
Furthermore, in the limit $t\to\infty$,
\[
\frac{X_t}{t^\half} \to N(0, {\cal D})
\]
in distribution.

\section{Without vertical advection}
\label{mix}

Throughout this section, we assume that the vertical advection $v=0$.
Furthermore, suppose that the horizontal advection has the form
\[
u(t,y) = u_1(t,y) + u_2(t,y),
\]
where $u_1$ and $u_2$ are both periodic functions of time with period
$2\pi/\omega$, and, for all $t$ and $y$, that $u'_1(t,y) = u'_1(t,a-y)$,
and $u'_2(t,y) = - u'_2(t,a-y)$, where again the dashes denote that the
$y$ means have been removed.  In this case, since \eqref{PDE} is
linear, $f$ can be split into $f_1$ and $f_2$ forced respectively by
$u'_1$ and $u'_2$, and inheriting the same symmetries in $y$.  Thus
\eqref{Taylor} implies that 
\[\label{Dmix}
{\cal D} = {\cal D}_1 + {\cal D}_2,
\]
where ${\cal D}_1$ and ${\cal D}_2$ are the Taylor dispersivities
respectively for $u_1$ and $u_2$, since $q=1/a$ and
\[
\int_0^a \frac{\partial f_1}{\partial y}\frac{\partial f_2}{\partial
  y} dy = 0,
\]
as the integrand is an odd function of $y - a/2$.

Note that \eqref{Dmix} does not hold in general for arbitrary $v$,
unless the periods of the two oscillations are different.

\section{Special cases}
\label{examples}

We now consider several example flows in which it is possible to find
the Taylor dispersivity exactly.  We begin by considering the way in
which shear- and pressure-driven flows interact in the case of a
Newtonian fluid, and then consider several cases of power-law fluid
flows.  In these examples, the vertical advection velocity $v$ is zero.
The former, in effect, provide solutions to the Hele-Shaw cell model
of the corresponding porous material system, namely a microscopic flow
driven by shearing the whole porous material, and a pressure-driven
flow for a fixed porous material. For biomechanical applications, see
Schmidt et. al. (2005).

In general we need to solve \eqref{PDE}, and the algebra is very messy
even in simple cases.  If $v$, the vertical advection velocity, is non-zero, even
if it is constant, the corresponding expressions for Taylor
dispersivity ${\cal D}$ are too complicated to show here, even for simple
horizontal flows.  However, using symbolic algebra packages, it is
possible to find simple exact expressions for ${\cal D}$ in some cases.

For simplicity in the rest of this section, we choose units that make
$\sigma = 1$, and $a=1$.  Also in these special cases $v = 0$, so
\eqref{PDE} reduces to
\[\label{PDE2}
\half \frac{\partial^2 f}{\partial y^2} + \frac{\partial f}{\partial
  t} = u'(t,y),
\]
with $\frac{\partial f}{\partial y}(t,0) = \frac{\partial f}{\partial
  y}(t,1) = 0$ for all $t$. Since $v=0$ and $a=1$, we find $q=1$,
which slightly simplifies \eqref{Taylor}.  With this
non-dimensionalization, the time-scale for the cross-channel
relaxation to occur is $1$, and the time-scale of the oscillations is
$\omega^{-1}$.  In the following examples, $u$ is proportional to
$\cos(\omega t)$. Thus, when doing the computer algebra, we can
consider $u$ and $f$ to be the real parts of functions proportional to
$\exp(i\omega t)$.

\begin{figure}
  \centering
  \includegraphics[angle=-90,width=\figurewidth]{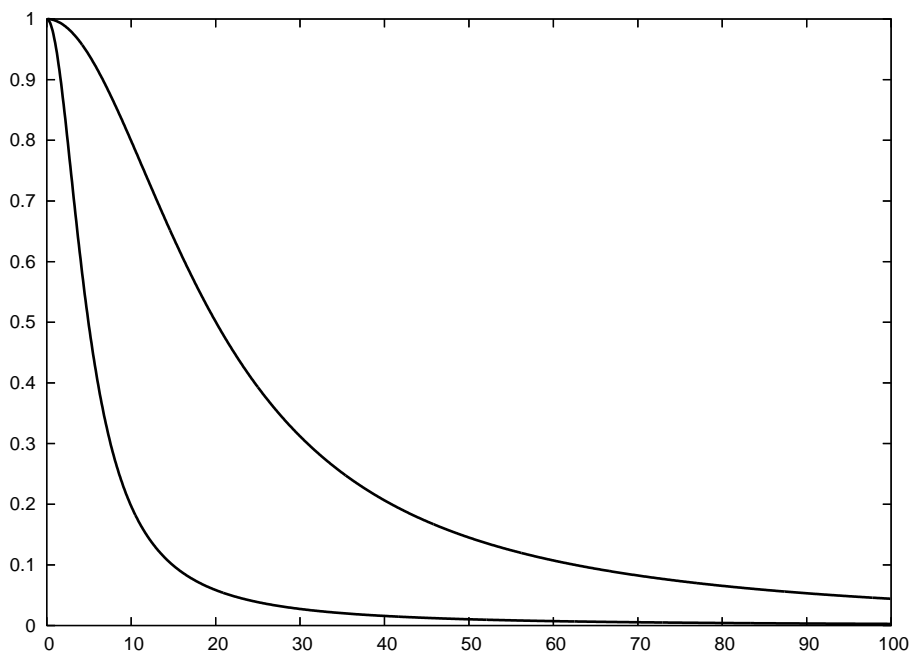}
  \caption{Plots of ${\cal D}(\omega)/{\cal D}(0)$. The upper curve is
    for $u(t,y) = \cos(\omega t) y(1-y)$ and the lower curve is for
    $u(t,y) = \cos(\omega t) y$.}
  \label{plot1}
\end{figure}

\subsection{$u(t,y) = \cos(\omega t) y$ and $v = 0$}

For $u(t,y) = \cos(\omega t) y$, using computer algebra, we find from
\eqref{PDE2} and \eqref{Taylor}
\[
\label{d1}
{\cal D} = d_1 \equiv \frac{\nu\cos(\nu) + \nu\cosh(\nu) - \sin(\nu) - \sinh(\nu)}
{2 \nu^5 (\cos(\nu) + \cosh(\nu))},
\]
where $\nu\equiv\omega^{\frac12}$ (see Figure \ref{plot1}).
In the limit of small $\omega$,
\[
{\cal D} \sim \frac1{60} - \frac{31\omega^2}{45360} +
\frac{5461\omega^4}{194594400} + \cdots.
\]
In the limit of large $\omega$,
\[
{\cal D} \sim \frac1{2\omega^2} + \cdots.
\]
This is in agreement with the result of Claes and Van de Broeck
(1991).  This is almost, but not exactly, the result previously found
by Young et. al. (1982) and extended by Jansons and Rogers (1995) for
the case of unbounded flow.  It differs by a factor of $2$ from the
result of Young et. al. (1982), which is explained by Jansons and
Rogers (1995), and lacks the oscillatory term of Jansons and Rogers
(1995). The large $\omega$ result here is, as expected, the average of
the result for the unbounded case.

\subsection{$u(t,y) = \cos(\omega t) y(1-y)$ and $v = 0$}

Similarly for $u(t,y) = \cos(\omega t) y(1-y)$, and solving for ${\cal D}$ using
computer algebra, we find
\[
\label{d2}
{\cal D} = d_2 \equiv \frac{\nu\cos(\nu) - \nu\cosh(\nu) - 3\sin(\nu) +
  3\sinh(\nu)} {6\nu^5(\cos(\nu) - \cosh(\nu))},
\]
(see Figure \ref{plot1}).
In the limit of small $\omega$,
\[
{\cal D} \sim \frac1{3780} - \frac{\omega^2}{1496880} +
\frac{\omega^4}{583783200} + \cdots.
\]
In the limit of large $\omega$,
\[
{\cal D} \sim \frac1{6\omega^2} + \cdots.
\]
This agrees with the result obtained by Claes and Van de Broeck (1991).

\subsection{A linear combination of the flows in the previous examples}

We now extend the work of Claes and Van de Broeck (1991), and consider
a flow driven by both a moving upper boundary and a pressure gradient,
but not necessarily in phase, namely
\[
u(t,y) = U_1\cos(\omega t) y + U_2 \cos(\omega t + \psi) y(1-y),
\]
where $\psi$ is a constant. From \eqref{Dmix}, \eqref{d1} and
\eqref{d2}, we find
\[
{\cal D} = U_1^2 d_1 + U_2^2 d_2.
\]
Notice that this result does not depend on the phase shift, which
could be important in some applications.

\subsection{Power-law fluids}

The same approach can be used to quickly find exact expressions for
some non-Newtonian flows. Here we consider
power-law fluid flows, driven by pressure gradients.  However, most of
the results are too complex to give here, and are essentially useful
only in computer programs.  However, with the procedure given
here, the reader can quickly recalculate them with the help of
computer algebra. These flows are of the following form:
\[
u_n(t,y) \equiv \cos(\omega t) (\textrm{const.} - |y - \half|^n),
\]
where the constant is fixed by the boundary conditions which may be
no-slip or partial slip. (Note that $n$ is not the defining parameter
of a power-law fluid, but it is more convenient here.) However the
results for the Taylor dispersivity does not depend on the constant.

Exact expressions for the Taylor dispersivity can be found for some
values of $n$.  However, the expressions for integer $n$ seem to be by
far the simplest, and we give a few examples of these below.  We were
unable to find the solution for general $n$.

As these flows are symmetrical about $y=\half$, we solve in a half range
which means we avoid the singularity at $y=\half$, which is problematic
for computer algebra.

In the following results, we do not take account of the relaxation
times for the fluid, so we must assume that the time-scale of the
oscillation is long compared with these.

\begin{figure}
  \centering
  \includegraphics[angle=-90,width=\figurewidth]{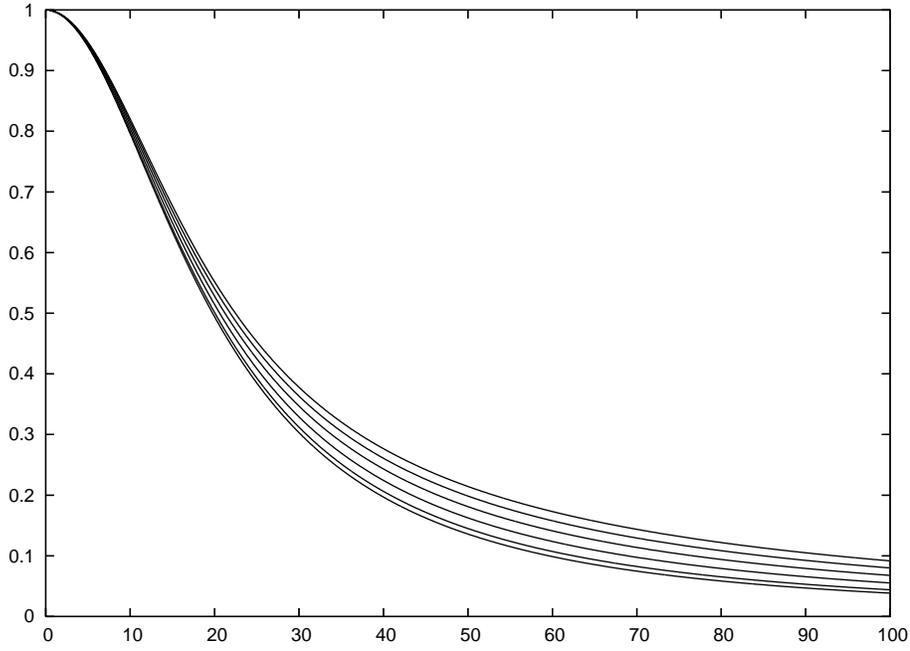}
  \caption{Plots of ${\cal D}(\omega)/{\cal D}(0)$ for $u_n$ with $n$
  increasing from $1$ at the bottom to
    $6$ at the top.}
  \label{plot2}
\end{figure}

The case $n=1$ is a rescaling of the Newtonian shear-flow result
above:
\[
{\cal D} = \frac{\nu\,\cos (\frac{\nu}{2}) + \nu\,\cosh (\frac{\nu}{2}) - 
     2\,\sin (\frac{\nu}{2}) - 2\,\sinh (\frac{\nu}{2})}{2\,\nu^5\,
     \left( \cos (\frac{\nu}{2}) + \cosh (\frac{\nu}{2}) \right) }.
\]

The case $n=2$ is for a Newtonian fluid, already considered above (see
\eqref{d2}):
\[
{\cal D} = d_2.
\]
 
The case $n=3$:
\[
\begin{aligned}
  {\cal D} &=
  (9\,( \nu\,( -80 + \nu^4 ) \,\cos (\nu) - 
  \nu\,( -80 + \nu^4 ) \,\cosh (\nu) \\
  &+ 
  80\,( -4 + \nu^2 ) \,\cosh (\frac{\nu}{2})\,
  \sin (\frac{\nu}{2}) - 5\,( -32 + 8\,\nu^2 + \nu^4 ) \,
  \sin (\nu) \\
  &+ 80\,( 4 + \nu^2 ) \,\cos (\frac{\nu}{2})\,
  \sinh (\frac{\nu}{2}) + 5\,( -32 - 8\,\nu^2 + \nu^4 )
  \,\sinh (\nu)
  )) \\
  &/(160\,\nu^9\,( \cos (\nu) - \cosh (\nu) )).
\end{aligned}
\]

The case $n=4$:
\[
\begin{aligned}
{\cal D} &= (\nu\,( -672 + \nu^4 ) \,\cos (\nu) - 
    \nu\,( -672 + \nu^4 ) \,\cosh (\nu) \\
&- 7\,( -144 + 24\,\nu^2 + \nu^4 ) \,\sin (\nu) + 
    7\,( -144 - 24\,\nu^2 + \nu^4 ) \,\sinh (\nu)) \\
  &/ (56\,\nu^9\, ( \cos (\nu) - \cosh (\nu) )).
\end{aligned}
\]

The case $n=5$:
\[
\begin{aligned}
{\cal D} &=  (5\,( \nu\,( 414720 - 13824\,\nu^4 + 5\,\nu^8) \,\cos (\nu) \\
&+ \nu\,( -414720 + 13824\,\nu^4 - 5\,\nu^8 ) \,\cosh (\nu) \\
&- 17280\,( -96 + 24\,\nu^2 + \nu^4 ) \,\cosh(\frac{\nu}{2})\,\sin
(\frac{\nu}{2}) \\
&- 45\,( 18432 - 4608\,\nu^2 - 768\,\nu^4 + 48\,\nu^6 + \nu^8 )\,\sin (\nu) \\
&+ 17280\,( -96 - 24\,\nu^2 + \nu^4 ) \,\cos(\frac{\nu}{2})\,\sinh
(\frac{\nu}{2}) \\
&+ 45\,( 18432 + 4608\,\nu^2 - 768\,\nu^4 - 48\,\nu^6 + \nu^8 ) \,\sinh (\nu) )) \\
&/ (4608\,\nu^{13}\,( \cos (\nu) - \cosh (\nu) ) ).
\end{aligned}
\]

The case $n=6$:
\[
\begin{aligned}
{\cal D} &= 
(9\,( \nu\,( 11827200 - 54560\,\nu^4 + 7\,\nu^8) \,\cos (\nu) \\
&+ \nu\,( -11827200 + 54560\,\nu^4 - 7\,\nu^8 )\,\cosh (\nu) \\
&- 77\,( 230400 - 38400\,\nu^2 - 2560\,\nu^4 + 80\,\nu^6 +\nu^8 )
\,\sin (\nu)\\
&+ 77\,( 230400 + 38400\,\nu^2 - 2560\,\nu^4 - 80\,\nu^6 + \nu^8)
\,\sinh (\nu) )) \\
&/ (39424\,\nu^{13}\,( \cos (\nu) - \cosh (\nu) )).
\end{aligned}
\]

Notice that in each case $\omega^2{\cal D}$ tends to a constant as
$\omega\to\infty$. 

There is little point in going beyond $n=6$ as it is almost plug flow,
with almost all the dispersion occurring near the boundaries.

The above results suggest that in the case of integer $n$ there is
probably some not too complex generating function for the Taylor
dispersivities, but \uppercase{Mathematica} was unable to work through
the necessary algebra.  However, since it would need computer algebra
to unpack the generating function, one is probably better off
computing the required Taylor dispersivity directly.

A plot of all the results of this section are given in Figure
\ref{plot2}.  Notice that even though we have considered only
integer $n$, the curves are sufficiently close that interpolation would
be expected to well approximate intermediate values over the
interesting part of the $\omega$ range.  Also, given that the case
$n=2$ is known, the fact that the other plots are close to this, and
occur in the expected order, is a good consistency check on the other
results.

\subsection{Comment on large and small $\omega$ limits}

From \eqref{PDE2} we can see that, for small $\omega$, the approximate
balance is
\[
\half \frac{\partial^2 f}{\partial y^2} = u'(t,y),
\]
and therefore the Taylor dispersivity is approximately the
time-average of the non-oscillatory Taylor dispersivities for the
instantaneous flows.  In the sinusoidal case this is just $\half$ of
the non-oscillatory Taylor dispersivity for the peak $u$.

For large $\omega$, the approximate balance is
\[
\frac{\partial f}{\partial t} = u'(t,y),
\]
giving an asymptotic Taylor dispersivity
\[
{\cal D} \sim \frac1{2\omega^2} \int_0^1
\left(\frac{dU}{dy}(y)\right)^2 dy, \qquad\mbox{as $\omega\to\infty$},
\]
where $U$ is defined by $u(t,y) = \cos(\omega t) U(y)$.

\section{Conclusions}

We have shown that the alternative method for determining the Taylor
dispersivity for an oscillatory flow lends itself to exact evaluation
using computer algebra.  Figure \ref{plot1} shows that for the
two simple flows considered, namely an oscillatory simple shear flow
and an oscillatory Poiseuille flow, both have a strong frequency-
dependent Taylor dispersivity.  These systems therefore could be used
for measuring particle diffusivities more accurately by choosing a
frequency that moves the corresponding Taylor dispersivities into the
sensitive range.  Alternatively, the same physical system could be
used for separating particles of different diffusivities.
Experimentally the oscillatory simple shear flow is particularly easy
to set up in a laboratory.

One interesting direction for future research is applications to
transport in bones, and other connective tissues, subjected to
oscillatory stresses (for example due to exercise or during ultra
sound therapy).  It appears that physical processes of the type
discussed here are essential for normal function, and a better
understanding could result in more effective treatment.  These
applications however involve much more complex geometries, although
would be expected to be qualitatively similar.  The solutions here are
in effect those for the Hele-Shaw cell model of a porous material.
Given we have shown that the Taylor dispersivities for shear- and
pressure-driven flows do not interact in a channel flow, it is a
natural question to ask if this carries over (at least approximately)
to these biomechanical applications, where both effects are likely to
be present.

It is clear that many more non-Newtonian examples are possible with
this approach.  However, it is less clear if there are any nice exact
expressions for flows with a non-zero vertical (i.e. cross-channel)
component. The intermediate expressions in these cases are much more
complex, with more branch cuts, but there does appear to be hope.
This is one area in which computer algebra is currently weak.

\end{document}